\documentclass{article}
\usepackage[]{algorithm2e}
\usepackage{amsmath}

\usepackage{graphicx}
\usepackage{amsfonts}
\usepackage{amssymb}
\usepackage{amsthm}
\usepackage{subfig}
\usepackage{rotating}
\usepackage{threeparttable}

\newlength{\tempheight}
\newlength{\tempwidth}
\newcommand{\rowname}[1]% #1 = text
{\rotatebox{90}{\makebox[\tempheight][c]{#1}}}

\newcommand{\where} {\noindent where~}

\newcommand{\columnname}[1]% #1 = text
{\makebox[\tempwidth][c]{#1}}

\newcommand{\ua}{\uparrow}
\newcommand{\nc}{\newcommand}
\nc{\da}{\downarrow} \nc{\hc}{\hat{c}} \nc{\hS}{\hat{S}}
\nc{\bra}{\langle} \nc{\ket}{\rangle} \nc{\eq}{equation (\ref}
\nc{\h}{\hat} \nc{\hT}{\h{T}}\nc{\be}{\begin{eqnarray}}
\nc{\ee}{\end{eqnarray}}\nc{\rd}{\textrm{d}}\nc{\e}{eqnarray}\nc{\hR}{\hat{R}}\nc{\Tr}{\mathrm{Tr}}
\nc{\tS}{\tilde{S}}\nc{\tr}{\mathrm{tr}}\nc{\8}{\infty}\nc{\lgs}{\bra\ua,\phi|}\nc{\rgs}{|\ua,\phi\ket}
\nc{\hU}{\hat{U}}\nc{\lfs}{\bra\phi|}\nc{\rfs}{|\phi\ket}\nc{\hZ}{\hat{Z}}\nc{\hd}{\hat{d}}\nc{\mD}{\mathcal{D}}
\nc{\bd}{\bar{d}}\nc{\bc}{\bar{c}}\nc{\mc}{\mathcal}\nc{\ea}{eqnarray}\nc{\mG}{\mathcal{G}}\nc{\bce}{\begin{center}}
\nc{\ece}{\end{center}}

\newcommand{\x}{x}

\renewcommand{\xi}{{\x^i}}

\newcommand{\Real}{\mathbb R}
\newcommand{\PP}{\mathbb P}

\newcommand{\R}{{\mathcal R}}

\begin{document}

\title{Quantifying Heteroskedasticity via Bhattacharyya Distance}
\author{M.~Hassan, M.~Hossny, D.~Creighton and S.~Nahavandi}

\maketitle

\begin{abstract}
Heteroskedasticity is a statistical anomaly that describes differing variances of error terms in a time series dataset. The presence of heteroskedasticity in data imposes serious challenges for forecasting models and many statistical tests are not valid in the presence of heteroskedasticity. Heteroskedasticity of the data affects the relation between the predictor variable and the outcome, which leads to false positive and false negative decisions in the hypothesis testing. Available approaches to study heteroskedasticity thus far adopt the strategy of accommodating heteroskedasticity in the time series and consider it an inevitable source of noise. In these existing approaches, two forecasting models are prepared for normal and heteroskedastic scenarios and a statistical test is to determine whether or not the data is heteroskedastic. 
\smallbreak
This work-in-progress research introduces a quantifying measurement for heteroskedasticity. The idea behind the proposed metric is the fact that a heteroskedastic time series features a uniformly distributed local variances. The proposed measurement is obtained by calculating the local variances using linear time invariant filters. A probability density function of the calculated local variances is then derived and compared to a uniform distribution of theoretical ultimate heteroskedasticity using statistical divergence measurements. The results demonstrated on synthetic datasets shows a strong correlation between the proposed metric and number of variances locally estimated in a heteroskedastic time series. 
%estimating variance trends, calculating the changes in variance slopes, and finally obtaining the average slope angles. Data is drawn from series of theoretical and real data sets. %The proposed measurement shows reliability for measuring and quantifying heteroskedasticity in comparison to the hypothesis and numerical tests of heteroskedasticity.
\end{abstract}

%\begin{keyword}
%Heteroskedasticity, Bhattacharyya distance, KL-Divergence
%%\texttt{elsarticle.cls}\sep \LaTeX\sep Elsevier \sep template
%%\MSC[2010] 00-01\sep  99-00
%\end{keyword}

%\end{frontmatter}

\section{Introduction}
Quantifying heteroskedasticity is a relatively new approach to study this statistical artefact. While heteroskedasticity is dealt with as an inevitable source of noise that must be accounted for in forecasting models, it becomes a noise source in signal processing and machine learning techniques \cite{Foi2009}. This new  uncontrollable source of noise then creates a new challenge to quantify heteroskedasticity. Consequently, early solutions for quantifying heteroskedasticity adopted one of two schools, change point detection and local parameter estimation. Change point detection methods \cite{Has13b} utilise the available heteroskedasticity tests to perform a binary segmentation of a heteroskedastic time series into smaller homoskedastic fragments.
\smallbreak 
Local parameter estimation methods utilises convolution with linear time invariant filters to estimate local variance at every sample based on its neighbours within a certain window $w$. There are two methods that adopt the local parameter estimation approach. Heteroskedasticity Variance Index (HVI) derived a variance of local variance as an indication of heteroskedasticity \cite{Has12}. Slope of Local Variance index (SoLVi) used the slope of the trend of estimated local variance derived by HVI as an indication of heteroskedasticity \cite{Has13a}. 
\smallbreak
An alternative approach to measure heteroskedasticity is to sample the estimated local variances in the time series. By doing this, a probability distribution $p_{\sigma^2}$ of the local variances can be derived. In theory, and as demonstrated in Figure~\ref{fig:heteropdf}, a homoskedastic time series should have a consistent local variance $\sigma^2$ over time. Consequently, the probability distribution of a homoskedastic time series should be unimodal and centred around $\sigma^2$. On the other hand, a heteroskedastic time series should, in theory, approach a uniform distribution covering a wide range of local variances. The ultimate heteroskedasticity time series should, in theory, feature a uniform distribution $\mathcal{U}\left(0, \infty\right)$. Therefore, measuring the distance between the probability distribution of the local variances $p_{\sigma^2}$ and the uniform distribution provides a quantified measure of heteroskedasticity. In this section, we propose heteroskedasticity measures based on probability distribution metrics. The heteroskedasticity quantified measure for a time series $y$ is defined as follows: %Bhattacharyya measures [?], KL divergence [?] and Mahalanobis [?] are presented. 
\begin{figure*}[t]
\includegraphics[width=.32\linewidth]{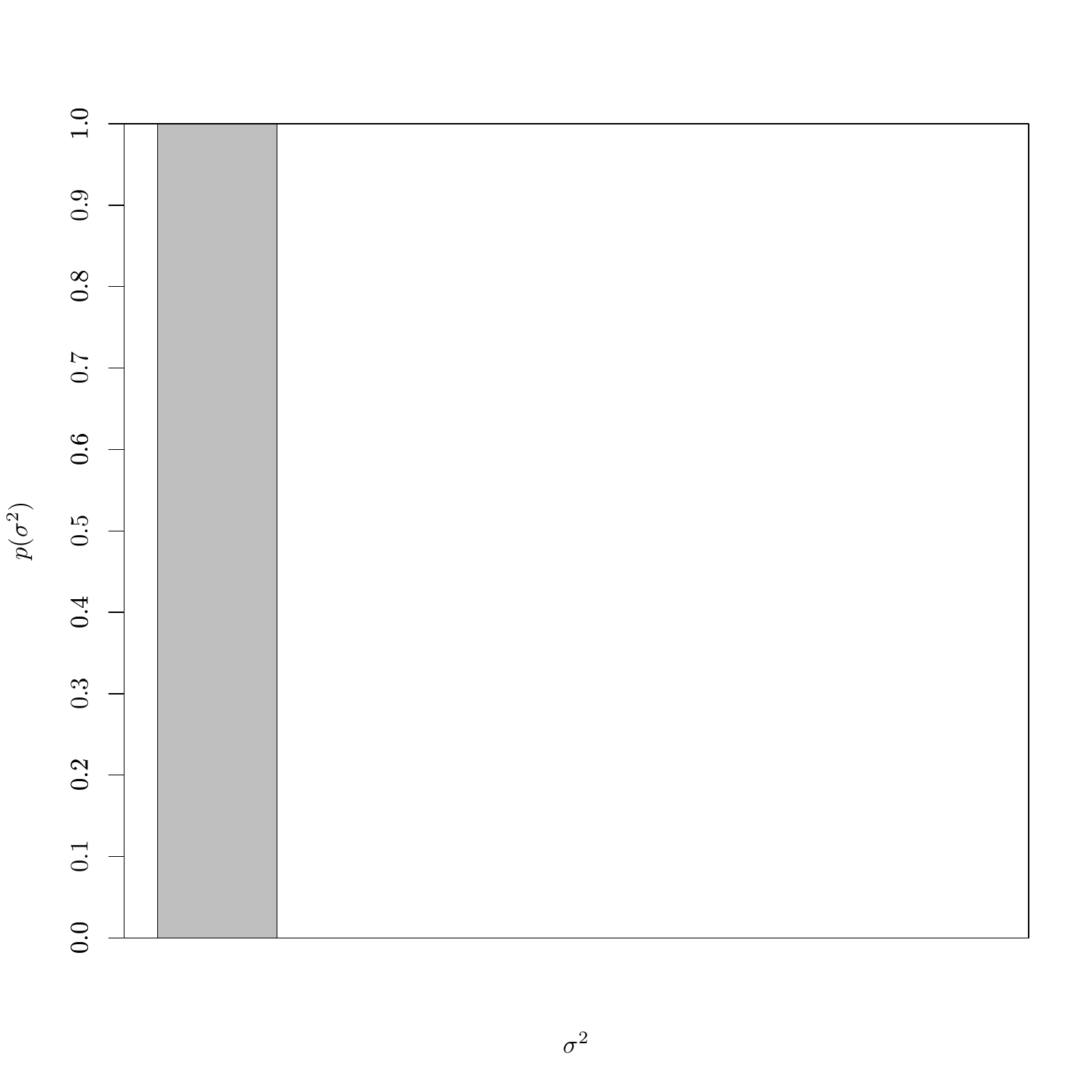}
\includegraphics[width=.32\linewidth]{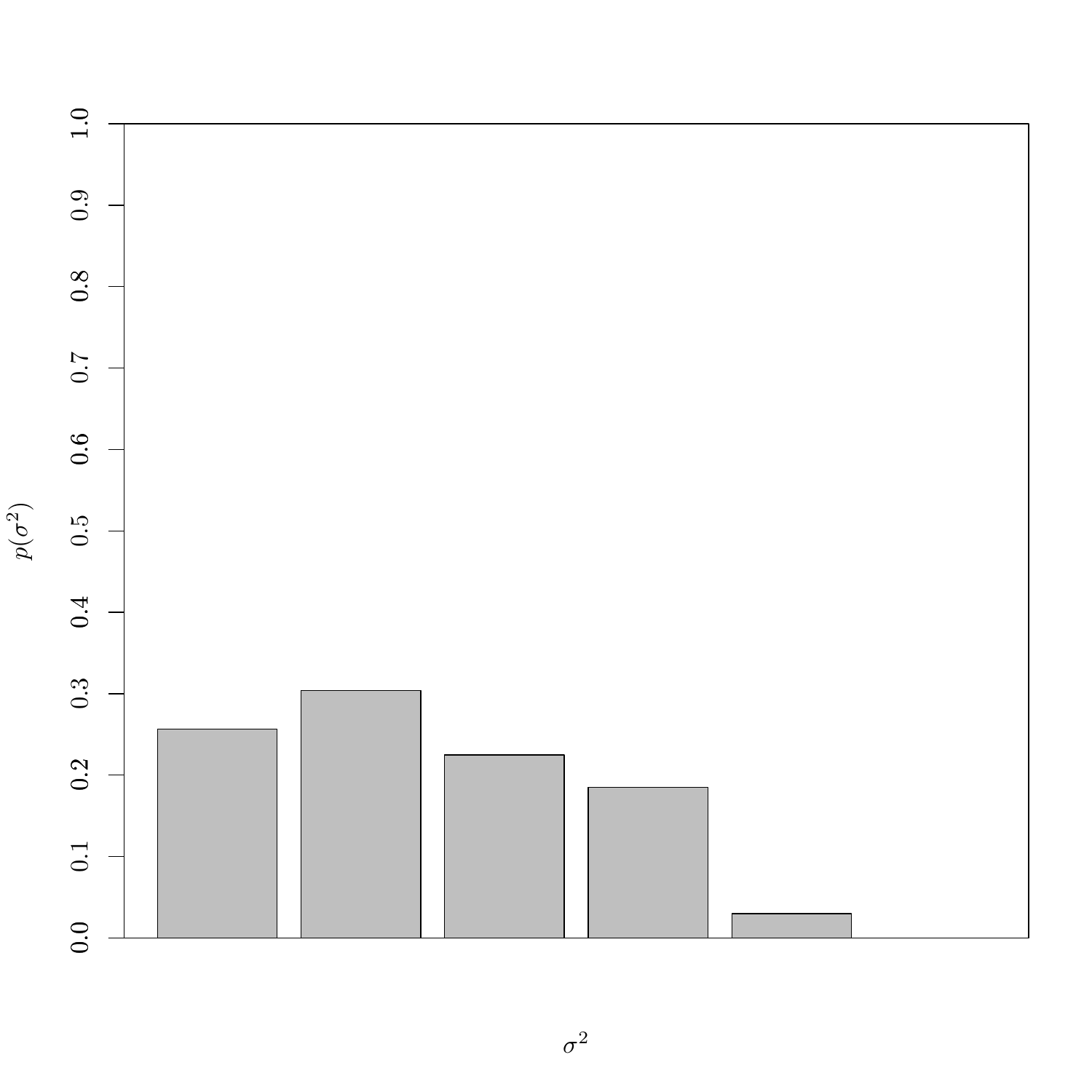}
\includegraphics[width=.32\linewidth]{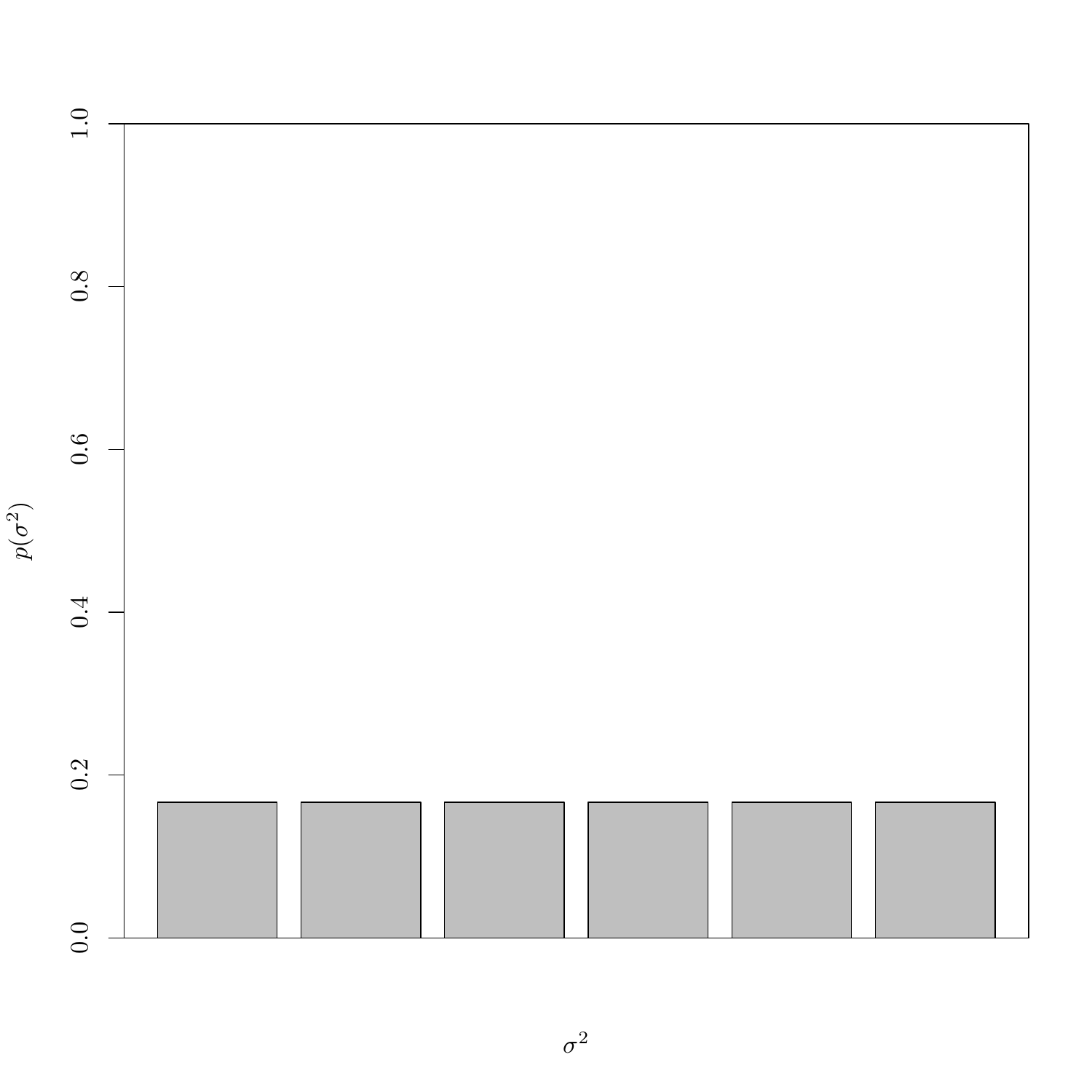}
\caption{Probability density function of local variances for a homoskedastic [left] and a heteroskedastic [middle] and a theoretically ultimate heteroskedastic [right] time series.}
\label{fig:heteropdf}
\end{figure*}

\begin{equation}
\mathcal{H}(y)=\Delta_p\left(P(\sigma^2_y), \mathcal{U}(0, \infty)\right),
\end{equation}

\where $\Delta_p:\PP^2 \to [0, 1]$ is a distribution distance function of the estimated local variances $\sigma^2_y$. Many probability distribution metrics are available. However, most of these metrics rely on entropies, joint probability density functions and sigma algebra. In this work-in-progress research, a justification for excluding three of the most famous probability distribution metrics such as mutual information, Tsallis mutual information, and Jensen-Shannon divergence are discussed. Additionally, a heteroskedasticity quantification function based on Bhattacharyya distance is introduced. 

\smallbreak
The rest of this paper is organised as follows. Section~\ref{sc:MI} covers mutual information and its variations. Section~\ref{sc:renyi} covers probability distribution divergence metrics based on Renyi entropy and introduces to the Bhattacharyya distance. Section~\ref{sc:bhatt} covers Bhattacharyya distance and the way it is utilised to quantify heteroskedasticity. Finally, Section~\ref{sc:conc} presents conclusion.

\section{Mutual Information (MI)}
\label{sc:MI}
Mutual Information between two random variables $X$ and $Y$ derives a cross entropy between the joint probability distribution $p(x, y)$ and the ultimate scenario of complete mutual independence $p(x) \cdot p(y)$ as follows:

\begin{equation}
\label{eq:MI}
MI(X;Y)=\int_{} \int_{} p(x, y)~\log\left(\frac{p(x, y)}{p(x)p(y)}\right) dx~dy.
\end{equation}

\noindent MI is used to measure the information shared between $X$ and $Y$ and equals to zero when $X$ and $Y$ are completely independent as follows:
 
\begin{eqnarray}
\label{eq:????}
MI(X;Y)&=&\int_{} \int_{} p(x) p(y)~\log\left(\frac{p(x)p(y)}{p(x)p(x)}\right)~dx~dy \nonumber \\
~~~~~~~&=&\int_{} \int_{} p(x) p(y)~\log 1~dx~dy=0. 
\end{eqnarray}

\noindent Mutual Information does not, however, provide a good solutions for quantifying heteroskedasticity because it is only bounded with the maximum entropy of $X$ or $Y$ as follows.

\begin{eqnarray}
\label{eq:????}
MI(X;X)&=&\int_{} p(x)~\log\left(\frac{p(x)}{p(x)p(x)}\right)~dx \nonumber \\
~~~~~~~&=&\int_{} p(x)~\log\left(\frac{1}{p(x)}\right)~dx \nonumber \\
%~~~~~~~&=&\int_{} p(x)~\log\left(\frac{1}{p(x)}\right)~dx\\
~~~~~~~&=&-\int_{} p(x)~\log p(x)~dx = H(X),
\end{eqnarray}

\where $H(X)$ is the entropy of $X$.

\subsection{Tsallis Driven Mutual Information (MI$_\alpha$)}
Another variation of mutual information was proposed by Cvejic \emph {et~al.} in \cite{Cve06}. They proposed to use the tunable Tsallis entropy \cite{Tsa88} described below.

\begin{equation}
\label{eq:TsallisMI}
MI_\alpha(X, Y)=\frac{1}{1-\alpha} \left(1-\int_{}\int_{} \frac{p(x)^\alpha}{p(y)^{1-\alpha}}~dx~dy \right),
\end{equation}

\where $\alpha \in \Real-\{1\}$ and $MI_\alpha(X, Y) \to MI(X, Y)$ as $\alpha \to 1$. This was proven by applying l'hopital rule on eq.~\ref{eq:TsallisMI} and substituting $\alpha=1$.

\begin{eqnarray}
\label{eq:????}
MI_1(X, Y)&=&\lim_{\alpha \to 1} MI_\alpha(X, Y) \nonumber \\
%~&=&\lim_{\alpha \to 1} \frac{1-\int\int p(x)^\alpha p(y)^{\alpha-1}~dx~dy}{1-\alpha} \\
~&=&\lim_{\alpha \to 1} \frac{\frac{d}{d\alpha}\left[1-\int\int p(x)^\alpha p(y)^{\alpha-1}~dx~dy\right]}{\frac{d}{d\alpha}\left[1-\alpha\right]}  \nonumber \\
%~&=&\lim_{\alpha \to 1} \int\int \frac{p(x)^\alpha p(y)^{1-\alpha} \ln p(x)}{p(x)^\alpha p(y)^{1-\alpha} \ln p(y)}~dx~dy\\
~&=&\lim_{\alpha \to 1} \int\int p(x)^\alpha p(y)^{1-\alpha}\frac{\ln p(x)}{\ln p(y)}~dx~dy  \nonumber \\
%~&=&\int\int p(x)\frac{\ln p(x)}{\ln p(y)}~dx~dy \nonumber \\
~&=&MI(X, Y).
\end{eqnarray}

\subsection{Jensen-Shannon Divergence}
Jensen-Shannon divergence metric uses sigma algebra \cite{Sha48} to derive an intermediate random variable $M=\frac{1}{2}(X+Y)$ which serves as a reference point to measure distance of $X$ and $Y$ from using mutual information as follows:

\begin{equation}
\label{eq:????}
JSD(X, Y)=\frac{1}{2}MI(X, M)+\frac{1}{2}MI(Y, M).
\end{equation}

\noindent While this metric is bounded to $0 \le JSD(X, Y) \le 1$, deriving the mixture distribution of the random variable $M$ is computationally intensive.

\section{Renyi Divergence}
\label{sc:renyi}
Renyi divergence \cite{Ren60} uses a generalised form of Shannon, Hartley, min-, and collision- entropies \cite{Sha48, Kon09} and is formulated as follows:

\begin{equation}
\label{eq:????}
H_\alpha=\frac{1}{1-\alpha} \log \left(\int_{}p(x)^\alpha\right).
\end{equation}

Renyi's divergence metric is then formulated as follows.

\begin{equation}
\label{eq:R_alpha}
\R_\alpha(X, Y)=\frac{1}{1-\alpha}\log \left(\int_{}\int_{} p(x)^\alpha p(y)^{1-\alpha}~dx~dy\right). 
\end{equation}

As Reynyi entropy generalises many entropies, its divergence metric also generalises many divergence metrics. For example, when $\alpha \to 1$ the Renyi entropy converges to Shannon's entropy and the divergence metric converges to the Mutual Information metric by applying l'hopital rule as follows:

\begin{eqnarray}
\label{eq:????}
H_1(X)&=&\lim_{\alpha \to 1} \frac{\log\left(\int_{}p(x)^\alpha\right)}{1-\alpha} \nonumber \\
%~~~~~~~~~~&=&\lim_{\alpha \to 1} \frac{\frac{\int p(x)^\alpha \log p(x)}{\int p(x)^\alpha}}{-1} \nonumber \\
~~~~~~~~~~&=&\lim_{\alpha \to 1} -\frac{\int p(x)^\alpha \log p(x)}{\int p(x)^\alpha} \nonumber \\
~~~~~~~~~~&=&-\frac{1}{\int p(x)} \int p(x) \log p(x) \nonumber \\
~~~~~~~~~~&=&-\int p(x) \log p(x) = H(X)\\
\R_1 &=&\int_{} \int_{} p(x, y)~\log\left(\frac{p(x, y)}{p(x)p(y)}\right) dx~dy \nonumber \\
&=&MI(X, Y)
\end{eqnarray}

\noindent Additionally, Renyi divergence also correlates with Bhattacharyya coefficient when $\alpha=\frac{1}{2}$ as follows:

\begin{eqnarray}
\label{eq:????}
\R_\frac{1}{2}(X)&=&\frac{\log \left(\int_{}\int_{} p(x)^\frac{1}{2} p(y)^{1-\frac{1}{2}}~dx~dy\right)}{1-\frac{1}{2}} \nonumber \\
~&=&-2\log \left(\int_{}\int_{} \sqrt{p(x)p(y)}~dx~dy\right) \nonumber \\
~&=&-2 \log BC(X, Y) \nonumber \\
~&=&2 \Delta^B_p(X, Y).
\end{eqnarray}

\section{Bhattacharyya Heteroskedasticity Measure}
\label{sc:bhatt}
\subsection{Bhattacharyya Distance}
Bhattacharyya-based metrics rely on deriving the Bhattacharyya Coefficient (BC)
\cite{Bha43}. The BC measures the closeness between two probability distributions $p$ and $q$ by measuring how disjoint they are as follows:

\begin{equation}
BC(p, q)=\sum_{x \in X}{\sqrt{p(x)q(x)}}
\end{equation}

\begin{figure}[t]
\centering
\includegraphics[width=.95\linewidth]{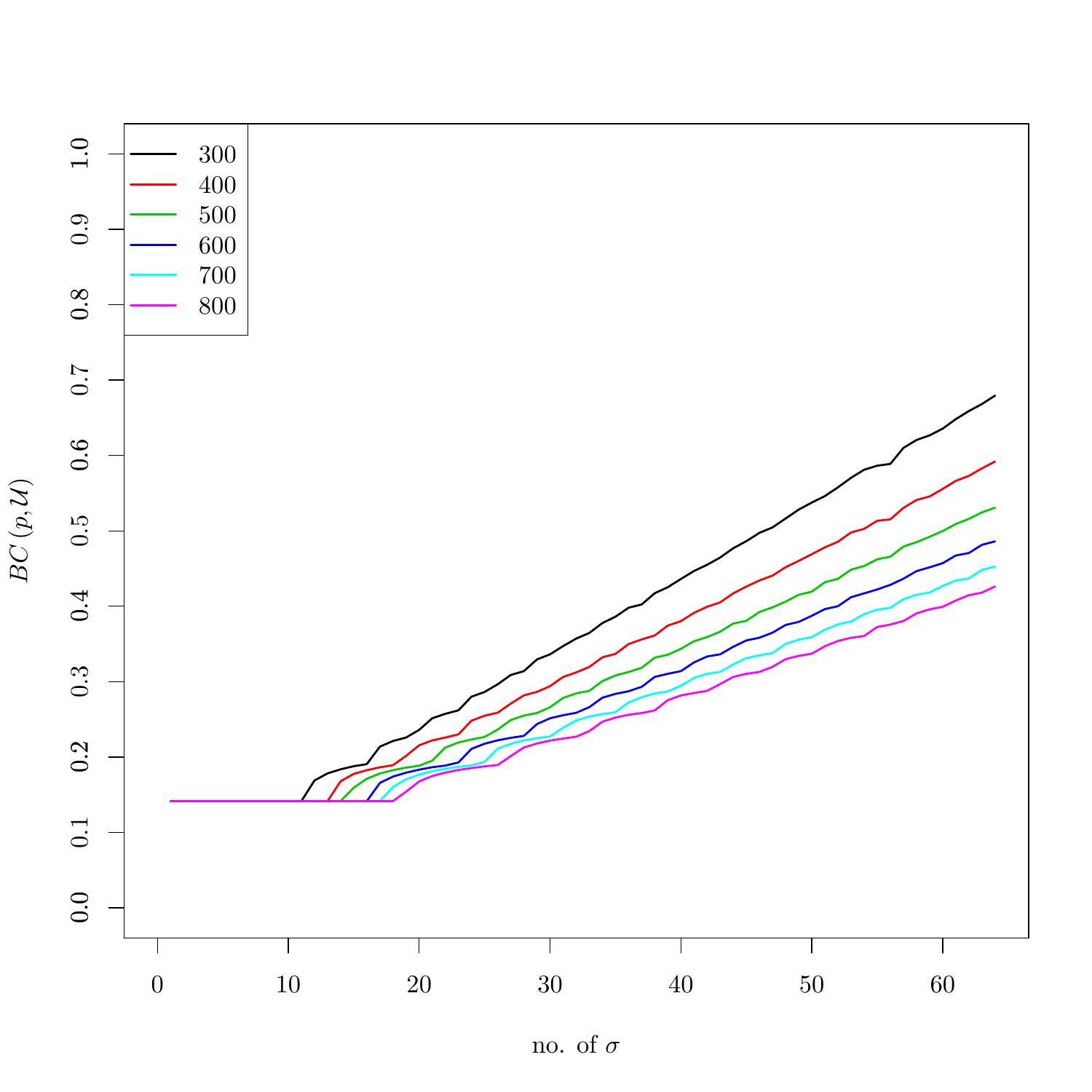}
\caption{Bhattacharyya coeffecient of a for time series generated using different 64 sigmas. The graphs demonstrate different kernel sizes $w$.}
\label{fig:BC}
\end{figure}

Figure~\ref{fig:BC} shows the Bhattacharyya coefficient as number of local variances increase in the dataset. Bhattacharyya coefficient has an upper bound of $1$ if and only if $p(x)=q(x)$. 

\begin{figure}[t]
\centering
\includegraphics[width=.95\linewidth]{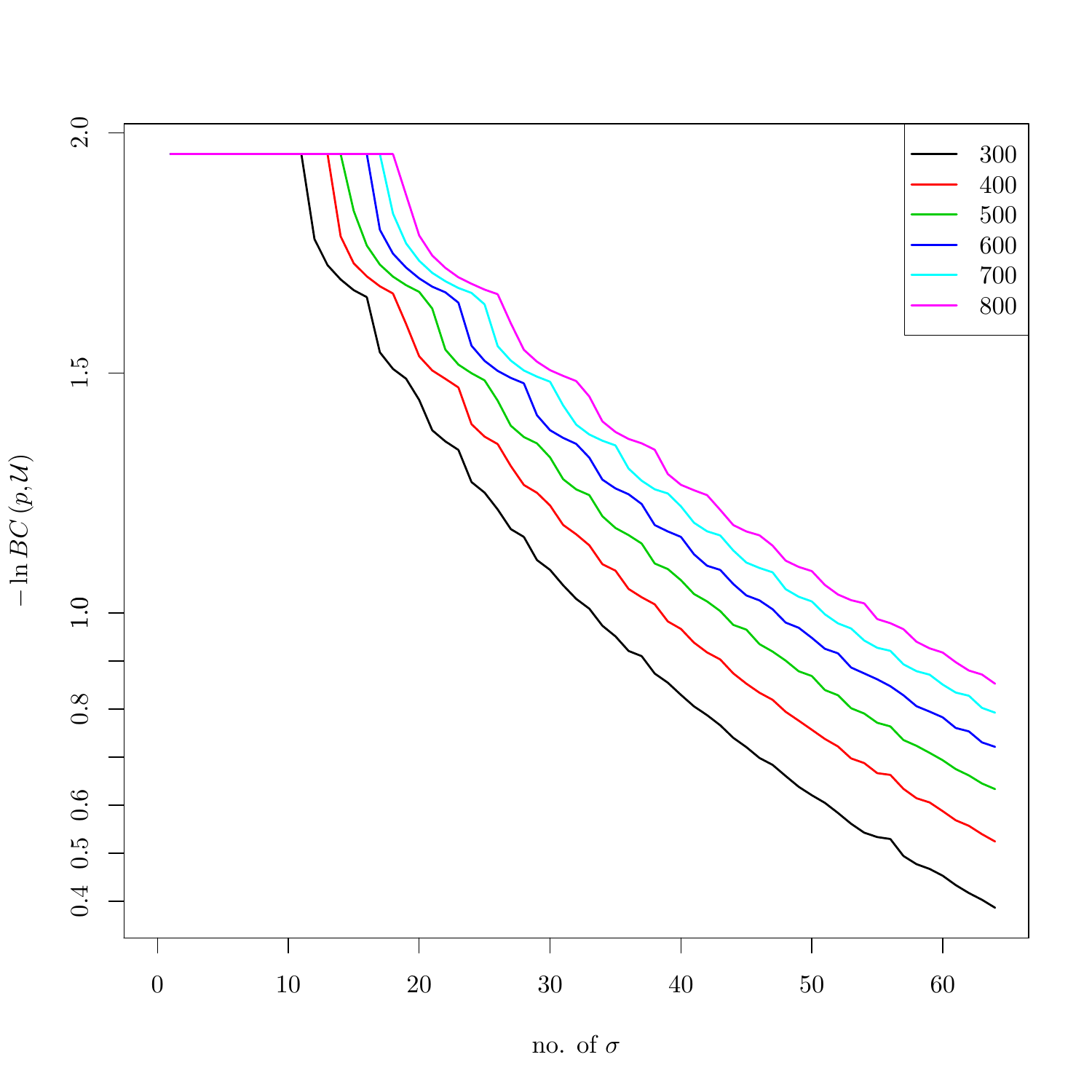}
\caption{Bhattacharyya distance of a for time series generated using different 64 sigmas. The graphs demonstrate different kernel sizes $w$.}
\label{fig:BD}
\end{figure}

This coefficient is then used to derive the Bhattacharyya distance as 

\begin{equation}
\Delta^B_p(p, q)=-\ln BC(p, q)
\end{equation}

\noindent However, this distance function has no upper bound and does not satisfy the triangulation inequality. Figure~\ref{fig:BD} demonstrates the Bhattacharyya distance. 

\subsection{Hellinger Distance}
Finally, Hellinger \emph {et~al.} provided a sound Bhattacharyya based divergence metric that is bounded and satisfies the triangulation inequality in \cite{Hel09}. The Hellinger metric is derived from Bhattacharyya coefficient as:

\begin{equation}
\Delta^H_p(p, q)=1-\sqrt{1-BC(p, q)}
\end{equation}

\begin{figure}[t]
\centering
\includegraphics[width=.95\linewidth]{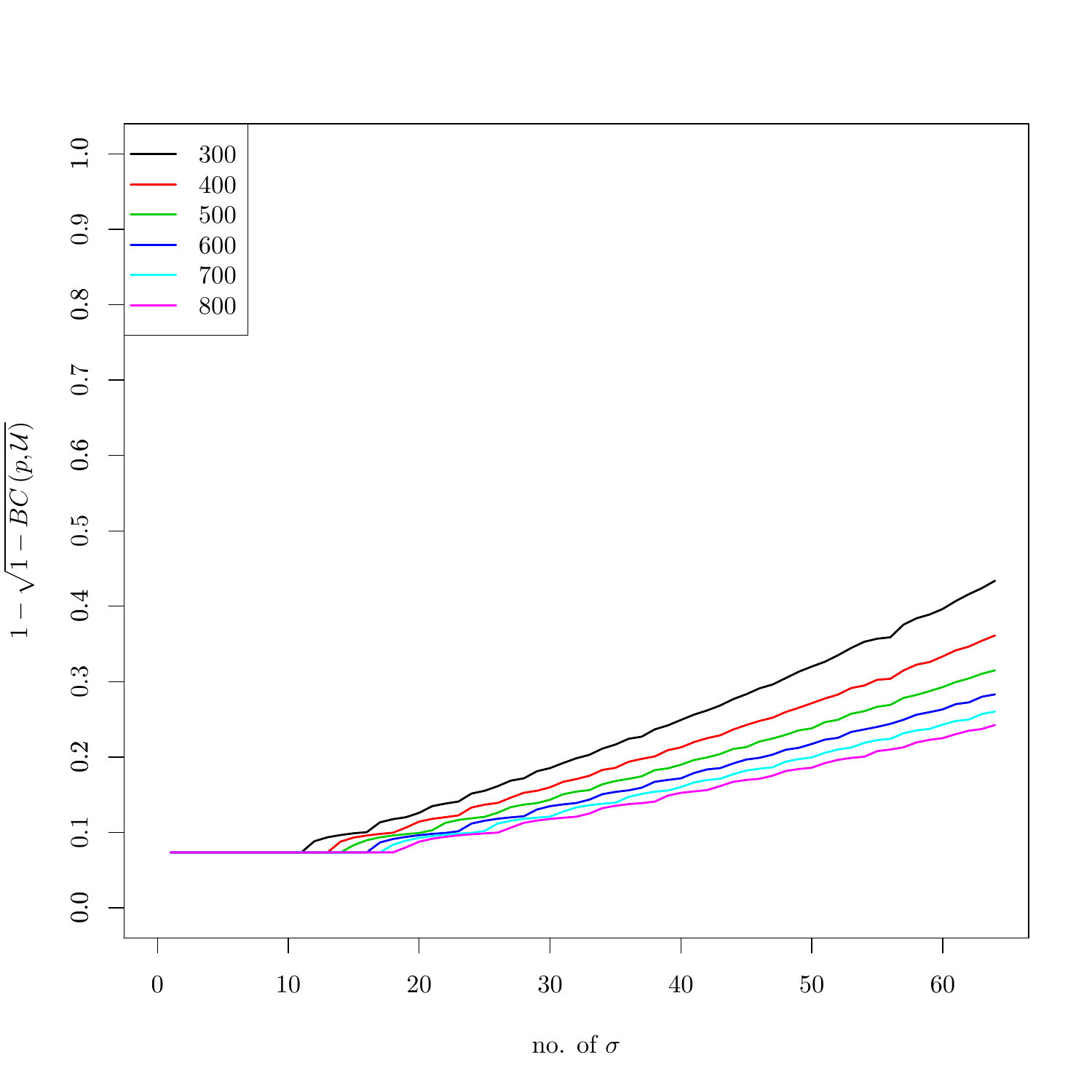}
\caption{Hellinger coefficient of a for time series generated using different 64 sigmas. The graphs demonstrate different kernel sizes $w$.}
\label{fig:Hellinger}
\end{figure}

Figure~\ref{fig:Hellinger} shows the effect of window size on Hellinger divergence metric.

\smallbreak
As a heteroskedastic time series, by definition, is derived from systems of different variances; the probability distribution of local variances $p(\sigma)$ of a heteroskedastic time series must be approaching a uniform distribution $\mathcal{U}$. On the other hand a homoskedastic time series will have a probability distribution further from the uniform distribution $\mathcal{U}$. To guarantee bounded function we chose Bhattacharayya coefficient over the Renyi driven metric in eq.~\ref{eq:R_alpha}. The Bhatacharayya heteroskedasticity measure is then formulated as follows:

\begin{equation}
\mathcal{H}_B(y)=\sum_{x \in X}\sqrt{P\left(\sigma^2_y\right) \mathcal{U}\left(\sigma^2_y\right)}
\end{equation}

\where $P\left(\sigma^2_y\right)$ is a probability distribution function of the estimated local variances $\sigma^2_y$. A Hellinger variation can also be derived with the same concept as follows:

\begin{equation}
\mathcal{H}_H(y)=1-\sqrt{1-\sum_{x \in X}\sqrt{P\left(\sigma^2_y\right) \mathcal{U}\left(\sigma^2_y\right)}}
\end{equation}

\section{Conclusions}
\label{sc:conc}
In this paper, we examine the divergence heteroskedasticity measures. Our motivation is that most of the available probability distribution metrics rely on entropies, joint density functions and sigma algebra. %Mutual information, Jensen-Shannon divergence and Renyi divergence were excluded. 
Measuring the distance between the probability distribution of the local variances $p_{\sigma^2}$ and the uniform distribution (ultimate heteroskedasticity) provides a quantified measure of heteroskedasticity. Consequently, the Bhattacharyya distance was adopted to introduce the Bhattacharyya heteroskedasticity measure. The main reason behind preferring the Bhattacharyya over the other KL-divergence measures is to guarantee a bounded function. The Bhattacharyya heteroskedasticity measure is then formulated using Hellinger variation to maintain the three propoerties of a distance function. %The proposed measure shows reliability in terms of measuring and quantifying heteroskedasticity.

\section*{Acknowledgement}
This research was fully supported by the Institute for Intelligent Systems Research and Innovation (IISRI).

%\section*{References}
%\bibliographystyle{els/model2-names}
%\printbibliography

\bibliography{QH15}

\end{document}